\documentclass[11pt,reqno]{amsart}
\usepackage{microtype}

\usepackage[dvipsnames]{xcolor}
\usepackage
[colorlinks=true,linkcolor=Maroon,citecolor=OliveGreen
]
{hyperref}
\usepackage{amsmath, amsfonts,amsthm,amssymb,verbatim,
    color,multirow,booktabs,mathdots,bm}
\usepackage[shortlabels]{enumitem}
\setlist[enumerate]{label={(\arabic*)}}
\usepackage[capitalize]{cleveref}
\crefname{equation}{}{}

\usepackage[abbrev, alphabetic]{amsrefs}  

\usepackage{soul}  

\newtheorem{theorem}{Theorem}

\newtheorem*{theorem*}{Theorem}

\theoremstyle{definition}

\theoremstyle{remark}

\def\Q{\mathbf{Q}}

\usepackage{todonotes}

\def\half{{\tfrac12}}

\numberwithin{equation}{section}

\author{Sean Eberhard}

\address{\parbox{\linewidth}{Mathematics Institute, Zeeman Building, University of Warwick, UK \vspace{0.1cm}}}
\email{sean.eberhard@warwick.ac.uk}

\thanks{SE is supported by the Royal Society.}

\begin{document}
\title{Ratios of consecutive values of the \mbox{divisor function}}

\begin{abstract}
    We show that the sequence of ratios $d(n+1) / d(n)$ of consecutive values of the divisor function attains every positive rational infinitely many times.
    This confirms a prediction of Erd\H{o}s.
\end{abstract}

\maketitle


Erd\H{o}s conjectured that the sequence of ratios $d(n+1) / d(n)$ of consecutive values of the divisor function is ``everywhere dense in $(0, \infty)$''~\cite{erdos86}.
In fact it follows easily from the generalized prime $k$-tuple conjecture that for every rational $x > 0$ there are infinitely many solutions to $d(n+1)/d(n) = x$.
We prove this unconditionally by building on the idea of Hasanalizade~\cite{hasanalizade} (see also Schlage-Puchta~\cite{schalge-puchta}), based on the following sieve theory result of Goldston, Graham, Pintz, and Y{\i}ld{\i}r{\i}m.

\begin{theorem}[\cite{GGPY}*{Corollary~2.1, special case $b_1=b_2=b_3=1$}]
    \label{thm:GGPY}
    Let $a_1, a_2, a_3, r_1, r_2, r_3$ be positive integers with $(r_i, a_i) = (r_i, a_i - a_j) = (r_i, r_j) = 1$ for all $i \ne j$.
    Let $L_i(x) = a_i x + 1$.
    Let $C$ be any positive integer.
    Then there are indices $i,j$ with $1 \le i < j \le 3$ such that there are infinitely many positive integers $x$ with both $L_i(x) / r_i, L_j(x) / r_j \in E_2(C)$,
    where $E_2(C)$ denote the set of products $p_1p_2$ where $p_1$ and $p_2$ are distinct primes and $p_1, p_2 > C$.
\end{theorem}

Set $(a_1, a_2, a_3) = (a, a+1, a+2)$, where $a$ is even. Then
\[
    a_2 L_1 - a_1 L_2 = a_3 L_2 - a_2 L_3 = \half a_3 L_1 - \half a_1 L_3 = 1.
\]
Let $r_1, r_2, r_3$ be pairwise coprime odd positive integers with $(r_i, a_i) = 1$ for each $i = 1, 2, 3$ and let $C$ be the largest prime factor of $a_1a_2a_3r_1r_2r_3$.
If $x$ is an integer such that $L_1(x)/r_1, L_2(x)/r_2 \in E_2(C)$ then $n = a_1 L_2(x)$ satisfies
\[
    \frac{d(n+1)}{d(n)} = \frac{d(a_2 L_1(x))}{d(a_1 L_2(x))} = \frac{d(a_2 r_1)}{d(a_1 r_2)}.
\]
Similarly, if $L_2(x) / r_2, L_3(x) / r_3 \in E_2(C)$ then $n = a_2 L_3(x)$ satisfies
\[
    \frac{d(n+1)}{d(n)} = \frac{d(a_3 L_2(x))}{d(a_2 L_3(x))} = \frac{d(a_3 r_2)}{d(a_2 r_3)},
\]
while if $L_1(x) / r_1, L_3(x) / r_3 \in E_2(C)$ then $n = \half a_1 L_3(x)$ satisfies
\[
    \frac{d(n+1)}{d(n)} = \frac{d(\half a_3 L_1(x))}{d(\half a_1 L_3(x))} = \frac{d(\half a_3 r_1)}{d(\half a_1 r_3)}.
\]
Thus, if $R$ denotes the set of values attained infinitely many times by the sequence $d(n+1)/d(n)$,
then \Cref{thm:GGPY} implies that $R$ contains at least one of the three values
\[
    \frac{d(a_2 r_1)}{d(a_1 r_2)},\qquad
    \frac{d(a_3 r_2)}{d(a_2 r_3)},\qquad
    \frac{d(\half a_3 r_1)}{d(\half a_1 r_3)}.
\]
Let $p_1, p_2, p_3$ be distinct primes such that $(a_1a_2a_3r_1r_2r_3, p_1p_2p_3) = 1$ and apply the fact just stated with $r'_i = r_i p_i^{e_i-1}$ in place of $r_i$, where $e_1, e_2, e_3 > 0$ are arbitrary positive integers. We deduce that $R$ contains one of
\[
    \frac{d(a_2 r_1) e_1}{d(a_1 r_2) e_2},\qquad
    \frac{d(a_3 r_2) e_2}{d(a_2 r_3) e_3},\qquad
    \frac{d(\half a_3 r_1) e_1}{d(\half a_1 r_3) e_3}.
\]
We can make the three expression above equal by setting
\begin{align*}
    e_1 &= d(a_2 r_1) d(a_3 r_2) d(\half a_1 r_3)^2, \\
    e_2 &= d(a_2 r_1) d(a_2 r_3) d(\half a_3 r_1) d(\half a_1 r_3), \\
    e_3 &= d(a_1 r_2) d(a_2 r_3) d(\half a_3 r_1)^2.
\end{align*}
It follows that
\[
    \frac{d(a_2 r_1) d(a_3 r_2) d(\half a_1 r_3)}
    {d(a_1 r_2) d(a_2 r_3) d(\half a_3 r_1)}
    \in R.
\]

Let $p_1, \dots, p_k, q_1, \dots, q_l$ be distinct odd primes. Let $x_i, y_i > 0$ ($1 \le i \le k)$ and $u_i, v_i > 0$ ($1 \le i \le l$) be positive integers and let
\begin{align*}
    a = 4 p_1^{x_1} \cdots p_k^{x_k} q_1^{u_1} \cdots q_l^{u_l},\qquad
    (r_1, r_2, r_3) = (1, p_1^{y_1} \cdots p_k^{y_k}, q_1^{v_1} \cdots q_l^{v_l}).
\end{align*}
Then $(r_i, a_i) = (r_i, 2) = (r_i, r_j) = 1$ for $i \ne j$, as required.
Thus $R$ contains
\[
    \frac{d(a_2 r_1) d(a_3 r_2) d(\half a_1 r_3)}
    {d(a_1 r_2) d(a_2 r_3) d(\half a_3 r_1)}
    =
    \frac43 \prod_{i=1}^k \frac{(x_i+1)(y_i+1)}{x_i+y_i+1} \prod_{i=1}^l \frac{u_i+v_i+1}{(u_i+1)(v_i+1)}.
\]
Let $f(x,y) = (x+1)(y+1) / (x+y+1)$. Noting that $4/3 = f(1, 1)$, we deduce that $R$ contains the subgroup $G$ of $\Q^\times_{>0}$ generated by $\{f(x,y) : x,y\ge 1\}$.

We claim that $G = \Q^\times_{>0}$. It suffices to prove that $G$ contains every prime number $p$, which we prove by induction on $p$. First, $2 = f(2,3) \in G$. Next, if $p = 2x+1 > 2$ then, since $x+1$ is a product of primes smaller than $p$, we have $p = (x+1)^2 / f(x, x) \in G$ by induction. This completes the proof.

\bibliography{refs}

\begin{bibdiv}
\begin{biblist}

\bib{erdos86}{inproceedings}{
      author={Erd\H{o}s, P.},
       title={Some problems on number theory},
        date={1986},
   booktitle={Proceedings of the seventeenth {S}outheastern international
  conference on combinatorics, graph theory, and computing ({B}oca {R}aton,
  {F}la., 1986)},
      volume={54},
       pages={225\ndash 244},
      review={\MR{885283}},
}

\bib{GGPY}{article}{
      author={Goldston, Daniel~A.},
      author={Graham, Sidney~W.},
      author={Pintz, Janos},
      author={Y{\i}ld{\i}r{\i}m, Cem~Y.},
       title={Small gaps between almost primes, the parity problem, and some
  conjectures of {E}rd{\H{o}}s on consecutive integers},
        date={2011},
        ISSN={1073-7928},
     journal={Int. Math. Res. Not. IMRN},
      number={7},
       pages={1439\ndash 1450},
         url={https://doi.org/10.1093/imrn/rnq124},
      review={\MR{2806510}},
}

\bib{hasanalizade}{article}{
      author={Hasanalizade, Elchin},
       title={Distribution of the divisor function at consecutive integers},
        date={2021},
        ISSN={0004-9727},
     journal={Bull. Aust. Math. Soc.},
      volume={103},
      number={2},
       pages={244\ndash 247},
         url={https://doi.org/10.1017/S0004972720000738},
      review={\MR{4229492}},
}

\bib{schalge-puchta}{article}{
      author={{Schlage-Puchta}, Jan-Christoph},
       title={{On a Problem by Erd{\H{o}}s and Mirsky on the ratio of the
  number of divisors of consecutive integers}},
        date={2025-03},
     journal={arXiv e-prints},
       pages={arXiv:2504.11463},
      eprint={https://arxiv.org/abs/2504.11463},
}

\end{biblist}
\end{bibdiv}

\end{document}